\documentclass[11pt]{article}
\usepackage{amsfonts}
\usepackage{amssymb}
\usepackage{epsfig}
\usepackage{amsmath}

\newcommand{\ignore}[1]{}
\ignore{The text in between these parentheses is ignored by LaTeX
ignore}

\def\@begintheorem#1#2{\par\bgroup{\sc #1\ #2. }\it\ignorespaces}
\def\@opargbegintheorem#1#2#3{\par\bgroup{\sc #1\ #2\ (#3). } \it\ignorespaces}
\def\@endtheorem{\egroup}
\newtheorem{theorem}{Theorem}[section]
\newtheorem{corollary}[theorem]{Corollary}
\newtheorem{lemma}[theorem]{Lemma}
\newtheorem{proposition}[theorem]{Proposition}
\newtheorem{problem}[theorem]{Problem}
\newtheorem{example}[theorem]{Example}
\newtheorem{algorithm}[theorem]{Algorithm}
\newtheorem{definition}[theorem]{Definition}
\newcommand{\bt}[1]{\begin{theorem}\label{#1}}
\newcommand{\bc}[1]{\begin{corollary}\label{#1}}
\newcommand{\bl}[1]{\begin{lemma}\label{#1}}
\newcommand{\bp}[1]{\begin{proposition}\label{#1}}
\newcommand{\bpro}[1]{\begin{problem}\label{#1}}
\newcommand{\be}[1]{\begin{example}\rm\label{#1}}
\newcommand{\ba}[1]{\begin{algorithm}\rm\label{#1}}
\newcommand{\bd}[1]{\begin{definition}\rm\label{#1}}

\newcommand{\et}{\end{theorem}}
\newcommand{\ec}{\end{corollary}}
\newcommand{\el}{\end{lemma}}
\newcommand{\ep}{\end{proposition}}
\newcommand{\epro}{\end{problem}}
\newcommand{\ee}{\end{example}}
\newcommand{\ea}{\end{algorithm}}
\newcommand{\ed}{\end{definition}}

\def\N{\mathbb{N}}

\def\R{\mathbb{R}}
\def\Z{\mathbb{Z}}

\def \G {{{\cal G}}}
\def \L {{{\cal L}}}

\def \IP {\mbox{IP}}

\newcommand{\boproof}{\noindent {\bf Proof. }}
\newcommand{\eoproof}{\hspace*{\fill} $\square$ \vspace{5pt}}
\newcommand{\red}{\sqsubseteq}
\newcommand{\FourTiTwo}{{\tt 4ti2}}

\begin{document}

\title{\bf N-Fold Integer Programming}
\author{
Jes\'us A. De Loera
\thanks{Supported in part by NSF grant DMS-0309694,
a 2003 UC-Davis Chancellor's fellow award, and
the Alexander von Humboldt foundation.}
\and
Raymond Hemmecke
\thanks{Supported in part by the European network ADONET 504438
and by grant FOR 468 of the Deutsche Forschungsgemeinschaft.}
\and
Shmuel Onn
\thanks{Supported in part by a grant from
ISF - the Israel Science Foundation, by the Technion President Fund,
and by the Jewish Communities of Germany Research Fund.}
\and
Robert Weismantel
\thanks{Supported in part by the European network ADONET 504438
and by grant FOR 468 of the Deutsche Forschungsgemeinschaft.}}

\date{}
\maketitle

\begin{abstract}
In this article we study a broad class of integer programming problems
in variable dimension. We show that these so-termed
{\em n-fold integer programming problems} are polynomial time solvable.
Our proof involves two heavy ingredients discovered recently:
the equivalence of linear optimization and so-called directed
augmentation, and the stabilization of certain Graver bases.

We discuss several applications of our algorithm to multiway
transportation problems and to packing problems. One important
consequence of our results is a polynomial time algorithm for the
$d$-dimensional integer transportation problem for long multiway tables.
Another interesting application is a new algorithm
for the classical cutting stock problem.
\end{abstract}

\section{Introduction}

The integer programming problem is the following discrete
optimization problem, where $\N$ denotes the set of nonnegative
integers, $A$ is an integer matrix and $b,c$ are integer vectors of
suitable dimensions:
$$
\min\, \{cx:\ Ax=b,\ x\in\N^q\}.
$$
It is well known to be generally NP-hard but polynomial time
solvable in fixed dimension $q$, see \cite{Sch}.

In this article, motivated by applications to high-dimensional
transportation problems and contingency tables and by the recently
discovered striking universality theorem for rational polytopes
\cite{DO2} (see Section \ref{Applications}), we study the following
class of integer programming problems in variable dimension.

\vskip.2cm\noindent
{\bf The n-fold integer programming problem}. Fix a $p\times q$
integer matrix $A$. Given positive integer $n$ and integer vectors
$b=(b^0,b^1,\dots,b^n)$ and $c=(c^1,\dots,c^n)$, where
$b^0\in\Z^q$, and $b^k\in \Z^p$ and $c^k\in\N^q$ for $k=1,\dots,n$,
find a nonnegative integer vector $x=(x^1,\dots,x^n)$, where $x^k\in\N^q$
for $k=1,\dots,n$, which minimizes $cx=\sum_{k=1}^n c^k x^k$ subject to
the equations $\sum_{k=1}^n x^k=b^0$ and $Ax^k=b^k$ for $k=1,\dots,n$.

\vskip.2cm\noindent
The term ``n-fold integer programming" refers to the problem being
almost separable into $n$ similar programs
$\min\{c^k x:\ Ax^k=b^k,\ x^k\in\N^q\}$ in fixed dimension; however, the
constraint $\sum_{k=1}^n x^k=b^0$ binds these programs together, and the
result is an integer program in large variable dimension $nq$.

\vskip.2cm\noindent
Let the {\em n-fold matrix} of $A$ be the following
$(q+np)\times nq$ matrix, with $I_q$ the $q\times q$ identity matrix:

$$A^{(n)}\quad:=\quad ({\bf 1}_n\otimes I_q)\oplus (I_n \otimes A)\quad=\quad
\left(
\begin{array}{ccccc}
  I_q    & I_q    & I_q    & \cdots & I_q    \\
  A  & 0      & 0      & \cdots & 0      \\
  0  & A      & 0      & \cdots & 0      \\
  \vdots & \vdots & \ddots & \vdots & \vdots \\
  0  & 0      & 0      & \cdots & A      \\
\end{array}
\right)\quad .
$$
\vskip.2cm\noindent
Then the n-fold integer programming problem can be
conveniently written in matrix form as
$$\min\{cx:\ A^{(n)}x=b,\ x\in\N^{nq}\}\quad.$$

In this article we establish the following theorem. Naturally,
the input size is $n$ plus the bit size of the integer objective vector
$c\in\Z^{nq}$ and the integer right-hand side vector $b\in\Z^{q+np}$.

\bt{Main}
Fix any integer matrix $A$. Then there is a polynomial time
algorithm that, given any $n$ and any integer vectors $b$ and $c$,
solves the corresponding n-fold integer programming problem.
\et

\noindent
The proof of this theorem involves two heavy ingredients. First, it makes
use of the equivalence of the linear optimization problem and the directed
augmentation problem, recently introduced and studied in \cite{SW}.
Second, it uses recent results of \cite{HS} and \cite{SS}
on the stabilization of certain Graver bases.

One important consequence of Theorem \ref{Main} is a polynomial time
algorithm for the $3$-way transportation problem for long tables,
settling its computational complexity;
see Section \ref{Applications} for details.

\vskip.2cm\noindent{\bf Corollary \ref{Threeway}\ }
Fix any $r,s$. Then there is a polynomial time algorithm that, given
$l$, integer objective vector $c$, and integer line-sums $(u_{i,j})$,
$(v_{i,k})$ and $(w_{j,k})$, solves the integer transportation problem
$$\min\{\,cx\ :\ x\in\N^{r\times s\times l}\,,\ \sum_i x_{i,j,k}=w_{j,k}
\,,\ \sum_j x_{i,j,k}=v_{i,k}\,,\ \sum_k x_{i,j,k}=u_{i,j}\,\}\ .$$

The n-fold integer programming problem
and theorem can be generalized as follows.

\vskip.2cm\noindent
{\bf Generalized n-fold integer programming}. Fix integer matrices $A,B$
of sizes $r\times q$ and $s\times q$, respectively. Given positive
integer $n$ and integer vectors $b=(b^0,b^1,\dots,b^n)$ and
$c=(c^1,\dots,c^n)$, with $b^0\in\Z^s$, and $b^k\in \Z^r$ and $c^k\in\N^q$
for $k=1,\dots,n$, find $x=(x^1,\dots,x^n)$ with $x^k\in\N^q$ for
$k=1,\dots,n$, which minimizes $cx=\sum_{k=1}^n c^k x^k$ subject to
$\sum_{k=1}^n Bx^k=b^0$ and $Ax^k=b^k$ for $k=1,\dots,n$.

\vskip.2cm\noindent
We have the following more general result, from which
Theorem \ref{Main} is deduced in the case $B=I_q$.

\vskip.2cm\noindent
\bt{Generalized} Fix any pair of integer matrices $A,B$ of compatible sizes.
Then there is a polynomial time algorithm that solves the generalized
n-fold integer programming problem on any input $n,b,c$.
\et

The article is organized as follows. In Section \ref{Applications} we discuss
applications of Theorem \ref{Main} to multiway transportation problems
and to some packing problems, as follows. In \ref{Applications1} we obtain
the aforementioned Corollary \ref{Threeway} which provides a polynomial time
solution to $3$-way integer transportation problems for long tables,
contrasting the recent universality theorem of \cite{DO2} for slim tables.
We also extend this result to $d$-way transportation problems for long
tables of any dimension (Corollary \ref{Multiway}). In \ref{Applications2}
we describe applications to a certain shipment problem
(Corollary \ref{Shipment}) and to the classical cutting stock problem
(Corollary \ref{CuttingStock}). In Sections 3 -- 4 we develop the
necessary ingredients for our n-fold integer programming algorithm,
as follows. In Section 3 we discuss Graver bases and augmentation.
We show (Lemma \ref{Graver}) that the Graver basis allows to solve the
directed augmentation problem introduced recently in \cite{SW}, and,
combining this with the results of \cite{SW}, show that any feasible solution
to an integer program can be augmented to an optimal one in polynomial time
provided the Graver basis is part of the input (Theorem \ref{GraverAug}).
In section 4 we discuss the stabilization of Graver bases discovered recently
in \cite{HS,SS}, and use it to show that Graver bases of n-fold matrices
can be computed in polynomial time (Theorem \ref{TheoremGraver}). Finally, in
Section 5, we combine all the ingredients from Sections 3 and 4, and prove our
main result Theorem \ref{Generalized} and its specialization Theorem \ref{Main}.

\section{Applications}
\label{Applications}

\subsection{High dimensional transportation problems}
\label{Applications1}

A {\em $d$-way transportation polytope} is the set of all
$m_1\times\cdots\times m_d$ nonnegative arrays $x=(x_{i_1,\dots,i_d})$
such that the sums of the entries over some of their lower dimensional
subarrays ({\em margins}) are specified. For simplicity of exposition, we shall
concentrate here only on $d$-way {\em line-sum polytopes}, of the form
$$T =  \left\{\,x\in\R_+^{m_1\times\cdots \times m_d}
\, : \, \sum_{i_1} x_{i_1,\dots,i_d}=u_{i_2,\dots, i_d}\,,
\ \sum_{i_2} x_{i_1,\dots,i_d}=u_{i_1,i_3,\dots, i_d}\,,
\ \dots\,,\ \sum_{i_d} x_{i_1,\dots,i_d}=u_{i_1,\dots,i_{d-1}}\,\right\}\,.$$
Transportation polytopes and their integer points (called
contingency tables by statisticians), have been studied and
used extensively in the operations research literature and in
the context of secure statistical data disclosure by public agencies
such as the census bureau and the national center for health statistics.
In the operations research literature, one is typically interested in the
integer and linear transportation problems, which are the integer and
linear programming problems over the transportation polytope, see
e.g. \cite{BR,KW,OR,QS,Vla,YKK} and references therein. In the statistics
community, one is often interested in the values an entry can attain
in all tables with fixed margins, related to the security of the entry
under margin disclosure, and in the construction of a Markov basis
allowing a random walk on the set of tables with fixed margins, see e.g.
\cite{AT,Cox1,Cox2,DO3,DFKPR,IJ} and references therein.

It is well known that the system defining a $2$-way transportation polytope
is totally unimodular. This implies that all the above problems are easy in
this case. However, already $3$-way transportation problems are much harder.
Consider the problem of deciding if a given $3$-way line-sum polytope of
$r\times s \times l$ arrays (with $r$ rows, $s$ columns and $l$ layers)
contains an integer point: the computational complexity of this problem
provides useful indication about the difficulty of the problems mentioned above.
If $r,s,l$ are all fixed, then the problem is solvable in
polynomial time by integer programming in fixed dimension $rsl$.
On the other hand, if $r,s,l$ are all variable part of the input,
then the problem is NP-complete \cite{IJ}. The in-between cases are
much more delicate. The case of two parameters $r,s$ variable and
one parameter $l$ fixed was recently resolved in \cite{DO1}, where it was
shown to be NP-complete, strengthening \cite{IJ}. Moreover, very recently,
in \cite{DO2}, the following striking universality result was shown.
\bp{Transportation}
Any rational polytope $P=\{y\in\R_+^n:Ay=b\}$ is polynomial time
representable as a $3$-way line-sum transportation polytope
of size $r\times s\times 3$ for some (polynomially bounded) $r$ and $s$,
$$T\quad=\quad\{\,x\in\R_+^{r\times s\times 3}\ :\ \sum_i x_{i,j,k}=w_{j,k}\,,\
\sum_j x_{i,j,k}=v_{i,k}\,,\ \sum_k x_{i,j,k}=u_{i,j}\,\}\ .$$
\ep
Here representable means that there is a coordinate-erasing projection from
$\R^{r\times s\times 3}$ onto $\R^n$ providing a bijection between $T$ and
$P$ and between the sets of integer points $T\cap\Z^{r\times s\times 3}$
and $P\cap\Z^n$. Thus, {\em any} rational polytope is an
$r\times s\times 3$ line-sum polytope, and {\em any} integer (respectively,
linear) programming problem is equivalent to an integer (respectively, linear)
$r\times s\times 3$ line-sum transportation problem. This result led to the
solution of several open problems from \cite{Vla} and \cite{YKK} and had
several implications on the complexity of Markov bases and the entry
security problem, see \cite{DO2} and \cite{DO3} for more details.

However, the last case, of two parameters $r,s$ fixed and one
parameter $l$ variable, has remained open and intriguing.
Here, as a consequence of Theorem \ref{Main}, we are able to
resolve this problem and show that both the decision and
optimization problems are polynomial time solvable.

\bc{Threeway}
Fix any $r,s$. Then there is a polynomial time algorithm that, given
$l$, integer objective vector $c$, and integer line-sums $(u_{i,j})$,
$(v_{i,k})$ and $(w_{j,k})$, solves the integer transportation problem
$$\min\{\,cx\ :\ x\in\N^{r\times s\times l}\,,\ \sum_i x_{i,j,k}=w_{j,k}
\,,\ \sum_j x_{i,j,k}=v_{i,k}\,,\ \sum_k x_{i,j,k}=u_{i,j}\,\}\ .$$
\ec
\boproof
We formulate the $3$-way integer transportation problem as an n-fold
integer program with $n:=l$, $p:=r+s$, and $q:=r\cdot s$.
Reindex the variables as $x^k_{i,j}:=x_{i,j,k}$ so that the variables
vector is $x=(x^1,\dots,x^n)$ with $x^k=(x^k_{i,j})\in\N^{r\times s}$
a $2$-way $r\times s$ table - the $k$th layer of the $3$-way table $x$.
Similarly write $c=(c^1,\dots,c^n)$ with $c^k\in\Z^{r\times s}$
for the objective vector. Next, put $b:=(b^0,b^1,\dots,b^n)$, with
$b^0\in\N^{rs}$ defined by $b^0:=(u_{i,j})$, and $b^k\in\N^{r+s}$
defined by $b^k:=((v_{i,k}),(w_{j,k}))$ for $k=1,\dots,n$.
Finally, let $A$ be the $p\times q=(r+s)\times r\cdot s$ matrix of
equations for the usual $2$-way transportation polytope, forcing
row-sums and column-sums on each of the $r\times s$ layers $x^k$
by $Ax^k=b^k$, $k=1,\dots,n$. Then the equations $Ax^k=b^k$ force the
line-sums $v_{i,k}$ and $w_{j,k}$, and the additional n-fold integer
program binding constraint $\sum_{k=1}^n x^k=b^0$
forces the ``long" line-sums $u_{i,j}$. This completes the encoding.
Since $r,s$ are fixed, so are $p,q$, and $A$, and therefore,
the corollary follows from Theorem \ref{Main}.
\eoproof
\be{Example1}
Consider the case $r=s=3$ (the smallest where the problem is
genuinely $3$-dimensional). Then $p=6$, $q=9$, and writing
$x^k=(x^k_{1,1},x^k_{1,2},x^k_{1,3},x^k_{2,1},x^k_{2,2},x^k_{2,3},
x^k_{3,1},x^k_{3,2},x^k_{3,3})$, the matrix $A$ which defines the
n-fold program providing the formulation of the
$3\times 3\times l$ transportation problem is
$$A\quad=\quad\left(
\begin{array}{ccccccccc}
  1 & 1 & 1 & 0 & 0 & 0 & 0 & 0 & 0 \\
  0 & 0 & 0 & 1 & 1 & 1 & 0 & 0 & 0 \\
  0 & 0 & 0 & 0 & 0 & 0 & 1 & 1 & 1 \\
  1 & 0 & 0 & 1 & 0 & 0 & 1 & 0 & 0 \\
  0 & 1 & 0 & 0 & 1 & 0 & 0 & 1 & 0 \\
  0 & 0 & 1 & 0 & 0 & 1 & 0 & 0 & 1 \\
\end{array}
\right)\quad.
$$
Already for this case, of $3\times 3\times l$ tables, the only polynomial
time algorithm for the corresponding line-sum integer transportation problem
we are aware of is the one guaranteed by Corollary \ref{Threeway} above.
\ee

Corollary \ref{Threeway} extends to transportation problems
of any dimension $d$, for {\em long} tables, namely, of size
$m_1\times\cdots\times m_{d-1}\times l$, where $m_1,\dots, m_{d-1}$
are fixed and only the length $l$ is variable, as follows.
\bc{Multiway}
Fix $d,m_1,\dots,m_{d-1}$. Then there is a polynomial time algorithm that,
given $l$, integer objective $c$, and line-sums
$(u_{i_2,\dots, i_d}),\dots, (u_{i_1,\dots, i_{d-1}})$,
solves the long multiway transportation problem
$$\min \{c x\ :\ x\in\N^{m_1\times\cdots \times m_{d-1}\times l}
\, : \, \sum_{i_1} x_{i_1,\dots,i_d}=u_{i_2,\dots, i_d}\,,
\ \dots\,,\ \sum_{i_d} x_{i_1,\dots,i_d}=u_{i_1,\dots,i_{d-1}}\,\}\,.$$
\ec
\boproof
The long multiway transportation problem can be encoded as an
n-fold integer program with $n:=l$, $p:=\sum_{i=1}^{d-1}m_i$, and
$q:=\prod_{i=1}^{d-1}m_i$, by reindexing the variables as
$x^{i_d}_{i_1,\dots,i_{d-1}}:=x_{i_1,\dots,i_d}$, letting $A$ be the
matrix of equations of line-sums of $(d-1)$-way transportation polytope
of $m_1\times\cdots\times m_{d-1}$ arrays, and proceeding in direct
analogy to the proof of Corollary \ref{Threeway}. The details are omitted.
\eoproof

\subsection{Some packing problems}
\label{Applications2}

\vskip.5cm\noindent{\sc Minimum cost shipment}

\vskip.2cm\noindent
The minimum cost shipment problem concerns the shipment of a large stock
of items of several types, using various vessels, with minimum possible cost.
More precisely, the data is as follows. There are $t$ types of items.
The weight of each item of type $j$ is $w_j$ and there are $n_j$ items
of type $j$ to be shipped. There are $v$ available vessels, where vessel
$k$ has maximum weight capacity $u_k$. The cost of shipping one item of
type $j$ on vessel $k$ is $p_{j,k}$. We now formulate this as an n-fold
integer programming problem. We set $n:=v$, $p:=1$, $q:=t+1$.
The defining matrix is the row vector $A=(A_j):=(w_1,w_2,\dots,w_t,1)$.
The variables vector is $x=(x^1,\dots, x^n)$
with $x^k=(x^k_1,\dots,x^k_t,x^k_q)$, where $x^k_j$ represents the number
of items of type $j$ to be shipped on vessel $k$ for $j=1,\dots,t$, and
$x^k_q$ is an extra slack variable representing the unused weight capacity
in vessel $k$. The cost vector is $c=(c^1,\dots, c^n)$ with
$c^k=(c^k_1,\dots,c^k_t,c^k_q)$, where $c^k_j:=p_{j,k}$ for $j=1,\dots,t$,
and $c^k_q:=0$. Finally, the demand vector is $b=(b^0,b^1,\dots,b^n)$ \break
with  $b^k:=u_k$ for $k=1,\dots,n$, and
$b^0:=(n_1,\dots,n_t,\sum_{k=1}^v u_k-\sum_{j=1}^t n_jw_j)$.
Then the resulting n-fold integer programming problem,
$\min\{cx:\ A^{(n)}x=b,\ x\in\N^{nq}\}$,
can be written in scalar form as follows:
\begin{eqnarray*}
\min & & \sum_{j=1}^q \sum_{k=1}^n c^k_j\,x^k_j
    \ = \  \sum_{j=1}^t \sum_{k=1}^v p_{j,k}\,x^k_j \\
s.t. & & \sum_{k=1}^n x^k_j=b^0_j=n_j\,,\quad\quad\ j=1,\dots, t \\
     & & \sum_{k=1}^n x^k_q=b^0_q=\sum_{k=1}^v u_k-\sum_{j=1}^t n_jw_j \\
     & & \sum_{j=1}^q A_j\,x^k_j = \sum_{j=1}^t w_j\,x^k_j + x^k_q
            = b^k = u_k \,,\quad\quad\ k=1,\dots, n\\
     & & x^k_j\in\N\,,\quad\quad\quad\quad j=1,\dots, q\,,\quad k=1,\dots, n\ .
\end{eqnarray*}
Assume that the number $t$ of types is fixed, but the numbers $n_j$ of items
of each type may be very large: this is a reasonable assumption in applications
(for instance, several types of automobiles to be shipped overseas, or several
types of appliances to be shipped on ground). Then we obtain the following
striking corollary of Theorem \ref{Main}, showing that the problem is
polynomial time solvable, where the input size is $v$ plus the bit size of
the integer numbers $n_j, u_k, p_{j,k}$ constituting the data.
Note that this result is {\em much stronger} than the standard
results on the {\em pseudo-polynomial time} solvability of this kind of
packing and knapsack-type problems using dynamic programming: our algorithm
can handle {\em very large} $n_j$ and $u_k$, possibly exponential
in the dimensional parameter $v$.

\bc{Shipment}
For any fixed number $t$ of types and type weights $w_j$, the minimum cost
shipment problem is solvable in time which is polynomial in the number $v$
of vessels and in the bit size of the integer numbers $n_j$ of items of each
type to be shipped, vessel capacities $u_k$, and shipment costs $p_{j,k}$.
\ec

\vskip.5cm\noindent{\sc The cutting stock problem}

\vskip.2cm\noindent
This is a classical manufacturing problem, where the usual setup is as follows:
a manufacturer supplies rolls of material (such as scotch-tape or band-aid)
in one of $t$ different widths $w_1,\dots,w_t$. The rolls are all cut out
from standard rolls of common large width $u$, coming out of the production line.
Given orders by customers for $n_j$ rolls of width $w_j$, the problem facing
the manufacturer is to meet the orders using the smallest possible number
of standard rolls. This is almost a direct special case of the minimum
cost shipment problem discussed above, with sufficiently many identical
vessels, say $v:=\sum_{j=1}^t\lceil n_j / \lfloor u / w_j \rfloor \rceil$,
of capacity $u_k:=u$ each, playing the role of the standard rolls, and with
cost $p_{j,k}:=w_j$ for each roll of width $w_j$ regardless of the standard
roll from which it is being cut out. The only correction needed is that
each slack variable $x^k_q=x^k_{t+1}$, measuring the unused width of
the $k$th standard roll, has cost of one unit instead of zero,
so that the total cost becomes the number of standard rolls used.
Thus the formulation as an n-fold program is with
$n:=\sum_{j=1}^t\lceil n_j / \lfloor u/w_j \rfloor \rceil$, $p:=1$, $q:=t+1$,
$A:=(w_1,w_2,\dots,w_t,1)$, variables $x^k_j$ representing the number
of rolls of width $w_j$ cut out of the $k$th roll for $j=1,\dots,t$
and $x^k_q$ representing the unused width of the $k$th standard roll,
costs $c^k_j:=w_j$ for $j=1,\dots,t$ and $c^k_q:=1$, and demands
$b^k:=u$ for $k=1,\dots,n$ and $b^0:=(n_1,\dots,n_t,nu-\sum_{j=1}^t n_jw_j)$.

Again, quite surprisingly, we get the following useful
corollary regarding this classical problem.

\bc{CuttingStock}
For any fixed $t$ and widths $w_1,\dots,w_t$,
the cutting stock problem is solvable in time polynomial in
$\sum_{j=1}^t\lceil n_j / \lfloor u / w_j \rfloor \rceil$ and in
the bit size of the numbers $n_j$ of orders and raw roll width $u$.
\ec

One common approach to the cutting stock problem makes use of so-called
{\em cutting patterns}, which are feasible solutions of the knapsack
problem $\{y\in\N^t\,:\,\sum_{j=1}^t w_jy_j\leq u\}$. This is useful
when the width $u$ of the standard rolls is of the same order of magnitude
as the demand widths $w_j$. However, when $u$ is much larger than the $w_j$,
the number of cutting patterns becomes prohibitively large to handle.
But then the values $\lfloor u/w_j\rfloor$ are large and hence
$n:=\sum_{j=1}^t\lceil n_j / \lfloor u/w_j \rfloor \rceil$
is small, in which case the result of Corollary \ref{CuttingStock}
using $n$-fold integer programming becomes particularly appealing.

\section{Graver bases and directed augmentation}

Consider the following family $\IP_A$ of integer programs in
standard form, with arbitrary demand vector $b\in\Z^m$ and
arbitrary objective vector $c\in\Z^n$, sharing the same constraint
matrix $A\in\Z^{m\times n}$,
\[
\IP_A(b,c):\quad \min\{cx:\ Ax=b,\ x\in\N^n\}.
\]
A {\em universal test set} for the family $\IP_A$ is a {\em finite} subset
$G$ of the lattice $\L(A):=\{x\in\Z^n:\ Ax=0\}$ of dependencies on $A$
such that whenever $x$ is
feasible but not optimal for $\IP_A(b,c)$ (where $b:=Ax$), there is an
{\em improving direction} $g\in G$, namely such that $x-g$ is feasible
and better, that is, $x-g\in\N^n$ and $cg>0$. Thus, a universal
test set enables the solution of the following augmentation problem.

\vskip.2cm\noindent{\bf Augmentation problem.}
Given $A\in\Z^{m\times n}$, $x\in\N^n$ and $c\in\Z^n$, either find an
improving direction $g\in\Z^n$, namely one with
$x-g\in \{y\in\N^n\,:\,Ay=Ax\}$ and $cg>0$, or assert that no such $g$ exists.

\vskip.2cm\noindent
An {\em augmentation oracle for a matrix $A$} is one that solves the
augmentation problem, that is, when queried on $x\in\N^n$ and $c\in\Z^n$,
it either returns an improving direction $g$ or asserts that none exists.
Clearly, an explicit universal test set $G$ for $A$ enables the efficient
realization of an augmentation oracle for $A$ by simply searching for
an improving direction $g\in G$. An oracle solving the augmentation problem,
and in particular, an explicit universal test set $G$, enable the following
simple iterative procedure that, for any program $\IP_A(b,c)$ in $\IP_A$ with
bounded objective function, converts any feasible $x$ to an optimal one:
``while there exists an improving direction $g$ set $x:=x-g$ and repeat''.

In 1975, Graver \cite{Graver:75} constructed, for \emph{every}
integer matrix $A$, a canonical universal test set. The {\em Graver basis}
${\cal G}(A)$ of $A$ can be defined as follows. First, we need to
extend the partial ordering $\leq$ from $\N^n$ to $\Z^n$. For
$u,v\in\Z^n$ we say that $u$ is {\em conformal} to $v$, denoted $u\red v$,
if $|u_i|\leq |v_i|$ and $u_iv_i\geq 0$ for $i=1,\ldots,n$,
that is, $u$ and $v$ lie in the same orthant of $\R^n$ and each component
of $u$ is bounded by the corresponding component of $v$ in absolute value.
With this, ${\cal G}(A)$ consists precisely of all $\red$-minimal vectors in
$\L(A)\setminus\{0\}$. For a more detailed introduction of Graver
bases we refer to \cite{Hemmecke:PSP}. The currently fastest algorithm
to compute Graver bases, based on a completion procedure and a
project-and-lift approach, is described in \cite{Hemmecke:SymmGraver}
and implemented in the software package \FourTiTwo\ \cite{4ti2}.

The Graver basis, being a universal test set, provides an augmentation oracle
and hence enables to convert any feasible solution to an optimal one for any
$\IP_A(b,c)$ by the iterative augmentation procedure above. But this in
itself is not enough to guarantee an {\em efficient} (polynomial time)
solution: a major remaining question is how many augmentation steps
are needed to reach an optimal solution.

Recently, in \cite{SW}, a directed version of the augmentation problem
was introduced; quite remarkably, it was shown that the number
of {\em directed} augmentation steps needed to reach optimality is polynomial.
We discuss this next. First, we define the directed augmentation problem.

\vskip.2cm\noindent{\bf Directed augmentation problem.}
Given $A\in\Z^{m\times n}$, $x\in\N^n$ and $c',c''\in\Z^n$, either
find $g=g^+-g^-\in\Z^n$ satisfying $x-g\in\{y\in\N^n\,:\,Ay=Ax\}$
and $c'g^+-c''g^->0$, or assert that none exists.

\vskip.2cm\noindent
Here and throughout, $g^+,g^-\in\N^n$ denote the {\em positive} and
{\em negative} parts of $g\in\Z^n$, defined by $g^+_i:=\max\{g_i,0\}$
and $g^-_i:=-\min\{g_i,0\}$ for $i=1,\dots,n$. Thus, the directed
augmentation problem involves {\em two} objective function vectors:
$c'$ controls the cost of the positive part of $g$ and $c''$
controls the cost of the negative part of $g$.
The usual augmentation problem occurs as the special case $c'=c''=c$.
A {\em directed augmentation oracle for $A$} is one that solves
the directed augmentation problem, i.e. when queried on $x\in\N^n$,
$c',c''\in\Z^n$, it either returns an improving direction $g$ or
asserts that none exists.

\vskip.2cm\noindent
In \cite{SW}, it was assumed that the input includes an upper bound vector
$u\in\N^n$ on the variables, so that the actual feasible set is
$\{y\in\N^n\,:\,Ay=Ax,\ y\leq u\}$. Under this assumption, the feasible
region is always bounded and there is always an optimal solution.
Further, the complexity estimates in \cite{SW} depended on the
bit size of $u$ and $b:=Ax$. However, this is not really needed.
Consider the integer program $\IP_A(b,c)\,:\ \min\{cy:\ Ay=b,\ y\in\N^n\}$
with $b=Ax$. Its objective function is bounded
(and hence there is an optimal solution) if and only if it is bounded for
the corresponding LP-relaxation $\min\{cy:\ Ay=b,\ y\in\R_+^n\}$,
which can be checked in polynomial time by linear programming.
Furthermore, whenever $\IP_A(b,c)$ has an optimal solution, it has one of bit
size polynomially bounded in the size of $A$ and $b=Ax$; this basically
follows from Cramer's rule, see e.g. \cite[Section 17.1]{Sch}.
Therefore, it is possible to compute an upper bound $u$ in terms of $A$
and $x$ only, and plug it into the analysis of \cite{SW}; for instance,
$u_k:=m!\,n(n+1)(\max|x_j|)(\max|A_{i,j}|)^m$ for $k=1,\dots,n$ will do.

With this, the results of \cite{SW} imply the following.
\bp{DAP}
There is a polynomial oracle time algorithm that, given $A\in\Z^{m\times n}$,
$x\in\N^n$, $c\in\Z^n$, solves the integer program $IP_A(b,c)$ with $b:=Ax$
by querying a directed augmentation oracle for $A$.
\ep
Here, as usual, {\em solving} the (feasible) integer program means that the
algorithm either returns an optimal solution or asserts that the objective
function is unbounded; and {\em polynomial oracle time} means that the number
of arithmetic operations, the number of calls to the oracle, and the size
of the numbers occurring throughout the algorithm are polynomially bounded
in the size of the input $A,x,c$.

It is not hard to see that if the matrix $A$ is totally unimodular, in
particular the incidence matrix of a directed graph, then a directed
augmentation oracle for $A$ can be realized using linear programming.
However, in general it is not clear which matrices $A$ admit
efficient directed augmentation.

As explained above, the Graver basis of a matrix $A$ yields an
augmentation oracle for $A$. We now show that, moreover, it enables the
realization of a directed augmentation oracle for $A$ as well.

\bl{Graver}
Let $\G(A)$ be the Graver basis of $A\in\Z^{m\times n}$. For any $x\in\N^n$ and
$c',c''\in\Z^n$, there is a $g\in\Z^n$ with $x-g\in\{y\in\N^n\,:\,Ay=Ax\}$
and $c'g^+-c''g^->0$ if and only if there is such $g\in \G(A)$.
\el

\boproof
Suppose $g$ is an improving direction. Then $g\in \L(A)\setminus\{0\}$
and hence can be written as a {\em conformal sum} of
(not necessarily distinct) elements of the Graver basis of $A$,
that is, $g=\sum g^i$ with $g^i\red g$ and $g^i\in\G(A)$ for all $i$.
To see this, recall that $\G(A)$ is the set of $\red$-minimal elements
in $\L(A)\setminus\{0\}$ and note that $\red$ is a well-ordering; if
$g\in\G(A)$, we are done; otherwise there is an $h\in\G(A)$ with
$h\sqsubset g$ in which case, by induction on $\red$, there is a
conformal sum $g-h=\sum g^i$ giving $g=h+\sum g^i$.

Now, $g^i\red g$ is equivalent to $(g^i)^+\leq g^+$ and $(g^i)^-\leq g^-$, so
the conformal sum $g=\sum g^i$ gives corresponding sums of the positive
and negative parts $g^+=\sum (g^i)^+$ and $g^-=\sum (g^i)^-$. Consequently,
\[
0\ <\ c'g^+ - c''g^- \ =\  c'\sum (g^i)^+ - c''\sum (g^i)^-
\ = \ \sum (c'(g^i)^+ - c''(g^i)^-)\ ,
\]
which implies that there is some $g^i$ in this sum with
$c'(g^i)^+ - c''(g^i)^->0$. Now, $g^i\in\G(A)\subset \L(A)$ so $Ag^i=0$
and hence $A(x-g^i)=Ax$. Finally we show that $x-g^i\geq 0$:
if $g^i_j\leq 0$ then $x_j-g^i_j\geq x_j\geq 0$; and if $g^i_j>0$ then
$g^i\red g$ implies $g^i_j\leq g_j$ and thus $x_j-g^i_j \geq x_j - g_j \geq 0$,
the last inequality holding because $g$ is an improving direction.
So $g^i\in\G(A)$ is an improving direction in the Graver basis.
\eoproof

As an immediate corollary of Proposition \ref{DAP}
and Lemma \ref{Graver}, we get the following statement.

\bt{GraverAug}
There is a polynomial time algorithm that, given any matrix $A\in\Z^{m\times n}$
along with its Graver basis $\G(A)$, and vectors $x\in\N^n$ and $c\in\Z^n$,
solves the integer program $IP_A(b,c)$ with $b:=Ax$.
\et

\noindent
While Theorem \ref{GraverAug} holds for any matrix, its complexity bound
depends on the size of the Graver basis which is part of the input.
Typically, the Graver basis is very large and its cardinality may be
exponential in $n$. However, in the next section we show that for a broad
and useful class of matrices, we can tame the behavior of the Graver basis,
leading to an efficient algorithm in terms of $A,x,c$ only.

\section{Graver bases of n-fold matrices}

Fix any pair of integer matrices $A$ and $B$ with the same number of columns,
of dimensions $r\times q$ and $s\times q$, respectively. The {\em n-fold
matrix of the ordered pair $A,B$} is the following $(s+nr)\times nq$ matrix,
\[
[A,B]^{(n)} \quad:=\quad ({\bf 1}_n\otimes B)\oplus (I_n\otimes A)\quad=\quad
\left(
\begin{array}{ccccc}
  B & B & B & \cdots & B      \\
  A & 0 & 0 & \cdots & 0      \\
  0 & A & 0 & \cdots & 0      \\
  \vdots & \vdots & \ddots & \vdots & \vdots \\
  0  & 0 & 0      & \cdots & A      \\
\end{array}
\right)\quad .
\]
With this, the generalized n-fold integer programming
problem can be conveniently written as
$$\min\{cx:\ [A,B]^{(n)}x=b,\ x\in\N^{nq}\}\quad.$$
The n-fold of a single matrix $A$, defined in the introduction, is the
special case $A^{(n)}=[A,I_q]^{(n)}$ with $B=I_q$ the $q\times q$ identity,
giving the regular (non-generalized) n-fold integer programming problem.

We now discuss a recent result of \cite{SS} and its extension in
\cite{HS} on the stabilization of Graver bases of n-fold matrices.
Consider vectors $x=(x^1,\ldots,x^n)$ with $x^k\in\N^q$ for $k=1,\dots,n$.
The {\em type} of $x$ is the number $|\{k\,:\,x^k\neq 0\}|$
of nonzero components $x^k\in\N^q$ of $x$. The following result of \cite{HS}
on the stabilization of Graver bases of $[A,B]^{(n)}$ extends
the earlier result for $B=I_q$ from \cite{SS}.

\bp{GraverComplexity}
For every pair of integer matrices $A\in\Z^{r\times q}$ and
$B\in\Z^{s\times q}$, there exists a constant $g(A,B)$
such that for all $n$, the Graver basis of $[A,B]^{(n)}$
consists of vectors of type at most $g(A,B)$.
\ep
The smallest constant $g(A,B)$ possible in the proposition
is called the {\em Graver complexity} of $A,B$.

Using Proposition \ref{GraverComplexity}, we now show that
$\G([A,B]^{(n)})$ can be computed in polynomial time.

\bt{TheoremGraver}
Fix any pair of integer matrices $A\in\Z^{r\times q}$ and $B\in\Z^{s\times q}$.
Then there is a polynomial time algorithm that, given $n$, computes the Graver
basis $\G([A,B]^{(n)})$ of the n-fold matrix $[A,B]^{(n)}$. In particular, the
cardinality and the bit size of $\G([A,B]^{(n)})$ are bounded by a polynomial
function of $n$.
\et
\boproof
Let $g:=g(A,B)$ be the Graver complexity of $A,B$ and consider any $n\geq g$.
We show that the Graver basis of $[A,B]^{(n)}$ is the union of
$n\choose g$ suitably embedded copies of the Graver basis of $[A,B]^{(g)}$.
Consider any $g$ indices $1\leq k_1<\dots <k_g\leq n$ and define a map
$\phi_{k_1,\dots,k_g}$ from $\N^{gq}$ to $\N^{nq}$ by sending
$x=(x^1,\ldots,x^g)$ to $y=(y^1,\ldots,y^n)$ defined by
$y^{k_t}:=x^t$ for $t=1,\dots,g$, and $y^k:=0$ for all other $k$.

We claim that the Graver basis of $[A,B]^{(n)}$ is the union of the images
of the Graver basis of $[A,B]^{(g)}$ under the $n\choose g$ maps
$\phi_{k_1,\dots,k_g}$ for all $1\leq k_1<\dots <k_g\leq n$, that is,
\begin{equation}
\G([A,B]^{(n)})\quad=\quad \bigcup_{1\leq k_1<\dots <k_g\leq n}
\phi_{k_1,\dots,k_g}(\G([A,B]^{(g)}))\quad .
\end{equation}
To see this, recall first that, by definition, the Graver basis of a
matrix $M$ is the set of all $\red$-minimal nonzero dependencies on $M$
(where a dependency on $M$ is a vector $x$ satisfying $Mx=0$).
Thus, if $x=(x^1,\ldots,x^g)\in\G([A,B]^{(g)})$ then $x$ is a
$\red$-minimal nonzero dependency on $[A,B]^{(g)}$, implying that
$\phi_{k_1,\dots,k_g}(x)$ is a $\red$-minimal nonzero dependency
on $[A,B]^{(n)}$ and hence $\phi_{k_1,\dots,k_g}(x)\in\G([A,B]^{(n)})$.
This establishes that the right-hand side of (1) is contained in the
left-hand side. Conversely, consider any $y\in\G([A,B]^{(n)})$. Then,
by Proposition \ref{GraverComplexity}, the type of $y$ is at most $g$,
so there are indices $1\leq k_1<\dots <k_g\leq n$ such that all
nonzero components of $y$ are among those of the reduced vector
$x:=(y^{k_1},\ldots,y^{k_g})$, and therefore $y=\phi_{k_1,\dots,k_g}(x)$.
Now, $y\in\G([A,B]^{(n)})$ implies that $y$ is a $\red$-minimal
nonzero dependency on $[A,B]^{(n)}$, and therefore $x$ is a $\red$-minimal
nonzero dependency on $[A,B]^{(g)}$ and hence $x\in\G([A,B]^{(g)})$,
showing that $y\in\phi_{k_1,\dots,k_g}(\G([A,B]^{(g)}))$. This establishes
that the left-hand side of (1) is contained in the right-hand side.
Thus, the Graver basis of $[A,B]^{(n)}$ is indeed given by (1).

Since $A,B$ are fixed and hence $g=g(A,B)$ is constant, the $g$-fold matrix
$[A,B]^{(g)}$ is also fixed and so the cardinality and bit size of its Graver
basis $\G([A,B]^{(g)})$ are constant as well. It follows from (1) that
$|\G([A,B]^{(n)})|\leq {n\choose g}|\G([A,B]^{(g)})|=O(n^g)$.
Further, each element of $\G([A,B]^{(n)})$ is an $nq$-dimensional vector
$\phi_{k_1,\dots,k_g}(x)$ obtained from some $x\in\G([A,B]^{(g)})$
(of constant bit size) by appending zero components, and therefore is
of linear bit size $O(n)$, showing that the bit size of the entire
Graver basis $\G([A,B]^{(n)})$ is $O(n^{g+1})$. Finally, it is clear that
the ${n\choose g}=O(n^g)$ images  $\phi_{k_1,\dots,k_g}(\G([A,B]^{(g)}))$
and their union $\G([A,B]^{(n)})$ can be computed in
time polynomial in $n$, completing the proof.
\eoproof

\be{Example2}
Consider the matrices $A=[1\ 1]$ and $B=I_2$.
The Graver complexity of the pair $A,B$ is $g(A,B)=2$. The 2-fold matrix
and its Graver basis, consisting of two antipodal vectors only, are
$$
[A,B]^{(2)}\ =\ A^{(2)}
\ =\ \left(
\begin{array}{cccc}
  1 & 0 & 1 & 0 \\
  0 & 1 & 0 & 1 \\
  1 & 1 & 0 & 0 \\
  0 & 0 & 1 & 1 \\
\end{array}
\right)\,,\quad\quad
\G([A,B]^{(2)})
\ =\ \pm \left(
\begin{array}{cccc}
  1 & -1 & -1 & 1 \\
\end{array}
\right)\quad.
$$
By Theorem \ref{TheoremGraver}, the Graver basis of the 4-fold matrix
$[A,B]^{(4)}=A^{(4)}$ can be computed by taking the union of the
images of the $6={4\choose 2}$ maps
$\phi_{k_1,k_2}:\N^{2\cdot 2}\longrightarrow\N^{4\cdot 2}$
for $1\leq k_1<k_2\leq 4$, and we obtain
$$
[A,B]^{(4)}=\left(
\begin{array}{cccccccc}
  1 & 0 & 1 & 0 & 1 & 0 & 1 & 0 \\
  0 & 1 & 0 & 1 & 0 & 1 & 0 & 1 \\
  1 & 1 & 0 & 0 & 0 & 0 & 0 & 0 \\
  0 & 0 & 1 & 1 & 0 & 0 & 0 & 0 \\
  0 & 0 & 0 & 0 & 1 & 1 & 0 & 0 \\
  0 & 0 & 0 & 0 & 0 & 0 & 1 & 1 \\
\end{array}
\right),\ \G([A,B]^{(4)})
=\pm \left(
\begin{array}{cccccccc}
  1 & -1 & -1 & 1 & 0 & 0 & 0 & 0 \\
  1 & -1 & 0 & 0 & -1 & 1 & 0 & 0 \\
  1 & -1 & 0 & 0 & 0 & 0 & -1 & 1 \\
  0 & 0 & 1 & -1 & -1 & 1 & 0 & 0 \\
  0 & 0 & 1 & -1 & 0 & 0 & -1 & 1 \\
  0 & 0 & 0 & 0 & 1 & -1 & -1 & 1 \\
\end{array}
\right)\, .
$$
\ee

\section{The polynomial time algorithm for n-fold integer programming}

We now provide the polynomial time algorithm for the generalized n-fold
integer programming problem
\begin{eqnarray}
\min\{cx:\ [A,B]^{(n)}x=b,\ x\in\N^{nq}\}\quad .
\end{eqnarray}
First, combining the results of the previous two sections, we get a
polynomial time procedure for converting any feasible solution to an
optimal one. We record this result in the following lemma.

\bl{Augmentation}
Fix any pair of integer matrices $A\in\Z^{r\times q}$ and $B\in\Z^{s\times q}$.
Then there is a polynomial time algorithm that, given $n$, objective vector
$c\in\N^{nq}$, and nonnegative integer vector $x\in\N^{nq}$, solves the
generalized n-fold integer programming problem in which $x$ is feasible,
i.e. the one with $b:=[A,B]^{(n)}x$.
\el
\boproof
First, apply the polynomial time algorithm underlying
Theorem \ref{TheoremGraver} on input $n$ and compute the Graver basis
$\G([A,B]^{(n)})$ of the n-fold matrix $[A,B]^{(n)}$. Then apply the
polynomial time algorithm underlying Theorem \ref{GraverAug} on input
$[A,B]^{(n)}$, $\G([A,B]^{(n)})$, $c$ and $x$, solving the integer program (2).
\eoproof

We now show that, moreover, given any $b$, we can efficiently
find an initial feasible solution to (2).

\bl{Feasible}
Fix any pair of integer matrices $A\in\Z^{r\times q}$ and $B\in\Z^{s\times q}$.
Then there is a polynomial time algorithm that, given $n$ and demand vector
$b\in\N^{s+nr}$, either finds a feasible solution $x\in\N^{nq}$ to the
generalized n-fold integer programming problem (2), or asserts that
no feasible solution exists.
\el
\boproof
Introduce $2n(r+s)$ auxiliary variables to the given generalized
n-fold integer program, and denote by $z$ the resulting vector
of $n(2(r+s)+q)$ variables. Consider the auxiliary integer program
of finding a nonnegative integer vector $z$ that minimizes the sum of
the auxiliary variables subject to the following system of equations,
with $I_r$ and $I_s$ the $r\times r$ and $s\times s$ identity matrices:
\[
\left(
\begin{array}{ccccccccccccccccc}
  B & I_s & -I_s &   0 &    0 & B & I_s & -I_s &   0 &    0 & B & \cdots & B & I_s & -I_s &   0 &    0\\
  A & 0 & 0 & I_r & -I_r & 0 & 0 & 0 & 0 & 0 & 0 & \cdots & 0 & 0 & 0 & 0 & 0\\
  0 & 0 & 0 & 0 & 0 & A & 0 & 0 & I_r & -I_r & 0 & \cdots & 0 & 0 & 0 & 0 & 0\\
  & & & & & & & & & & & \ddots & & & \\
  0 & 0 & 0 & 0 & 0 & 0 & 0 & 0 & 0 & 0 & 0 & \cdots & A & 0 & 0 & I_r & -I_r\\
\end{array}
\right)\cdot z\ =\ {\large b}\quad.
\]
This auxiliary program is in fact again a generalized n-fold integer program,
with matrices $\bar{A}=(A,0,0,I_r,-I_r)$ and $\bar{B}=(B,I_s,-I_s,0,0)$.
Since $A$ and $B$ are fixed, so are $\bar{A}$ and $\bar{B}$. Due to the
special structure of the auxiliary program, a feasible solution of this
program can be written down easily in terms of $b$. Consequently, the
auxiliary program can be solved by the algorithm underlying
Lemma \ref{Augmentation}, in time polynomial in $n$ and the bit size of $b$.
Since the auxiliary objective is bounded below by zero, the algorithm
will output an optimal solution $z$. If the optimal objective value is
(strictly) positive, then the original n-fold program (2) has no feasible
solution, whereas if the optimal value is zero, then the restriction of $z$ to
the original variables is a feasible solution $x$ of the original program (2).
\eoproof

Combining the results of Lemmas \ref{Augmentation} and \ref{Feasible},
we obtain the main result of this article.

\vskip.2cm\noindent
{\bf Theorem \ref{Generalized}\ }
{\it Fix any pair of integer matrices $A,B$ of compatible sizes.
Then there is a polynomial time algorithm that solves the generalized
n-fold integer programming problem on any input $n,b,c$.}

\vskip.2cm\noindent
Clearly, Theorem \ref{Main} is deduced from Theorem \ref{Generalized} in
the special case $B=I_q$. We emphasize again that, by solving the generalized
n-fold integer programming problem, we mean in the complete sense that
the algorithm concludes with precisely one of the following:
it either asserts that there is no feasible solution, or asserts that
the objective function is unbounded, or returns an optimal solution.

\noindent {\small Jesus De Loera}\newline
\emph{University of California at Davis, Davis, CA 95616, USA}\newline
\emph{email: deloera{\small @}math.ucdavis.edu},
\ \ \emph{http://www.math.ucdavis.edu/{\small $\sim$deloera}}

\noindent {\small Raymond Hemmecke}\newline
\emph{Otto-von-Guericke Universit\"at Magdeburg,
D-39106 Magdeburg, Germany}\newline
\emph{email: hemmecke{\small @}imo.math.uni-magdeburg.de},
\ \ \emph{http://www.math.uni-magdeburg.de/{\small $\sim$hemmecke}}

\noindent {\small Shmuel Onn}\newline
\emph{Technion - Israel Institute of Technology, 32000 Haifa, Israel}\newline
\emph{email: onn{\small @}ie.technion.ac.il},
\ \ \emph{http://ie.technion.ac.il/{\small $\sim$onn}}

\noindent {\small Robert Weismantel}\newline
\emph{Otto-von-Guericke Universit\"at Magdeburg,
D-39106 Magdeburg, Germany}\newline
\emph{email: weismantel{\small @}imo.math.uni-magdeburg.de},
\ \ \emph{http://www.math.uni-magdeburg.de/{\small $\sim$weismant}}

\end{document}